\documentclass{amsart}

\usepackage{txfonts}
\usepackage{amsmath,amssymb,amsfonts,enumerate,amsthm, amscd}

\renewcommand{\Im}{\mbox{Im}\,}

\newtheorem{theorem}{Theorem}[section]
\newtheorem{lemma}[theorem]{Lemma}
\newtheorem{proposition}[theorem]{Proposition}
\newtheorem{corollary}[theorem]{Corollary}

\theoremstyle{definition}
\newtheorem{definition}[theorem]{Definition}
\newtheorem{definitions}[theorem]{Definitions}
\theoremstyle{remark}
\newtheorem{remark}[theorem]{Remark}
\newtheorem{remarks}[theorem]{Remarks}
\newtheorem{example}[theorem]{Example}

\theoremstyle{Definition and Notation}

\begin{document}
\bibliographystyle{amsplain}

%\date{}
\title[]{On (strongly) Gorenstein (semi)hereditary rings}

\author{Najib Mahdou}
\address{Najib Mahdou\\Department of Mathematics, Faculty of Science and Technology of Fez, Box 2202, University S.M. Ben Abdellah Fez, Morocco.}

\author{Mohammed Tamekkante}
\address{Mohammed Tamekkante\\Department of Mathematics, Faculty of Science and Technology of Fez, Box 2202, University S.M. Ben Abdellah Fez, Morocco.}



\begin{abstract}

In this paper, we introduce and study the rings of Gorenstein
homological dimensions small or equal than  1, which we call
Gorenstein (semi)hereditary rings,  specially particular cases of
these rings, which we call strongly Gorenstein (semi)hereditary
rings.

\end{abstract}

\maketitle

%%%%%%%%%%%%%%%%%%%%%%%%%%%%%%%%%%%%%%%%%%%%%%%%%%%%%%%%%
%%%%%%%%%%%%%%%%%%%%%%%%%%%%%%%%%%%%%%%%%%%%%%%%%%%%%%%%%
%%%INTRODUCTION%%%%%%%%%%%%%%%%%%%%%%%%%%%%%%%%%%%%%%%%%%
\section{Introduction}
Throughout this paper, all rings are commutative with identity
element, and all modules are unital.\\

$\mathbf{\quad Setup\; and\; Notation:}$ Let $R$ be a ring, and
let $M$ be an $R$-module. As usual we use $pd_R(M)$, $id_R(M)$ and
$fd_R(M)$ to denote, respectively, the classical projective,
injective and flat dimensions of $M$. By $gldim(R)$ and $wdim(R)$
we denote, respectively, the classical global and weak dimensions
of R.\\
It is by now a well-established fact that even if R to be
non-Noetherian, there exists Gorenstein projective, injective and
flat dimensions of M, which are usually denoted by $Gpd_R(M)$,
 $Gid_R(M)$ and $Gfd_R(M)$, respectively. Some references are
 \cite{Bennis and Mahdou1, Bennis and Mahdou2, Christensen, Christensen
 Frankild, Enochs, Enochs2, Eno Jenda Torrecillas, Holm}.\\

 Recently in \cite{Bennis and Mahdou2}, the authors started the study of
global Gorenstein dimensions of rings, which are called, for a
commutative ring $R$, Gorenstein global  projective, injective,
and weak dimensions of $R$, denoted by $GPD(R)$, $GID(R)$, and
$G.wdim(R)$, respectively; and respectively, defined as
follows:\bigskip

$\begin{array}{cccc}
  1) & GPD(R) & = & sup\{ Gpd_R(M)\mid M$ $R-module\} \\
  2) & GID(R) & = & sup\{ Gid_R(M)\mid M$ $R-module\} \\
  3) & G.wdim(R) & = & sup\{ Gfd_R(M)\mid M$ $R-module\}
\end{array}$
\\

They proved that, for any ring R, $ G.wdim(R)\leq GID(R) = GPD(R)$
(\cite[Theorem 1.1 and Corollary 1.2(1)]{Bennis and Mahdou2}). So,
according to the terminology of the classical theory of
homological dimensions of rings, the common value of $GPD(R)$ and
$GID(R)$ is called
Gorenstein global dimension of $R$, and denoted by $G.gldim(R)$.\\
They also proved that the Gorenstein global and weak dimensions
are refinement of the classical global  and weak dimensions of
rings. That is : $G.gldim(R) \leq gldim(R)$ and $G.wdim(R)\leq
wdim(R)$ with equality if $wdim(R)$ is finite (\cite[Corollary
1.2(2 and 3)]{Bennis and Mahdou2}).\bigskip

In \cite{Bennis and Mahdou1}, the authors studied particular cases
of  Gorenstein projective, injective and flat modules which they
call strongly Gorenstein projective, injective and flat modules
respectively, and defined as follows:
\begin{definitions}\
\begin{enumerate}
    \item A module $M$ is said to be strongly Gorenstein projective ($SG$-projective for short), if
there exists an exact sequence of   projective modules  of the
form:
$$\mathbf{P}=\ \cdots\rightarrow P\stackrel{f}\rightarrow
P\stackrel{f} \rightarrow P\stackrel{f}
         \rightarrow P \rightarrow\cdots$$  such that  $M \cong \Im(f)$
and such that $Hom(-,P)$ leaves the sequence $\mathbf{P}$ exact whenever $P$ is projective.\\
  The exact sequence $\mathbf{P}$ is
called a strongly complete projective resolution.
    \item The strongly Gorenstein injective modules
are defined dually.
    \item A module $M$ is said to be strongly
Gorenstein flat ($SG$-flat for short), if there exists an exact
sequence of flat module  of the form:
$$\mathbf{F}=\ \cdots\rightarrow F\stackrel{f}\rightarrow F \stackrel{f}\rightarrow
F \stackrel{f}\rightarrow F \rightarrow\cdots$$ such that  $M
\cong \Im(f)$ and such that $-\otimes I$ leaves $\mathbf{F}$ exact
whenever $I$ is injective.

The exact sequence $\mathbf{F}$ is called a strongly complete flat
resolution.
\end{enumerate}
\end{definitions}
The principal role of the strongly Gorenstein projective and
injective modules is to give a simple characterization of
Gorenstein projective and injective modules, respectively, as
follows:
\begin{theorem}[\cite{Bennis and Mahdou1}, Theorem 2.7] A  module is Gorenstein projective (resp., injective)
if, and only if, it is a direct summand of a strongly Gorenstein
projective (resp., injective) module.

\end{theorem}
Using \cite[Theorem 3.5]{Bennis and Mahdou1} together with
\cite[Theorem 3.7]{Holm}, we have the next result:
\begin{proposition} Let $R$ be a coherent ring. A module is Gorenstein flat
if, and only if, it is a direct summand of a strongly Gorenstein
flat module.
\end{proposition}

 In this
paper, we often use the following Lemma:
\begin{lemma}\label{lemma1}
Consider the following diagram of modules aver a ring $R$.
$$(\star)\begin{array}{ccccccc}
  0\rightarrow & M & \stackrel{\alpha}\rightarrow & P & \stackrel{\beta}\rightarrow & M & \rightarrow 0 \\
   & _{u}\downarrow &  &  &  & _{u}\downarrow &  \\
  0\rightarrow & Q & \stackrel{\iota}\rightarrow & Q\oplus Q & \stackrel{j}\rightarrow & Q & \rightarrow 0 \\
\end{array}$$
where $M$ is a Gorenstein projective module, $P$ and $Q$ are
projective and $\iota$ and $j$ are the canonical injection and
projection respectively. Then, there is a morphism
$\gamma:P\rightarrow Q\oplus Q$ which complete $(\star)$ and make
it commutative.
\end{lemma}
\begin{proof}
If we apply the functor $Hom(-,Q)$ to the short exact sequence
$$(\ast)\quad 0\rightarrow  M  \stackrel{\alpha}\rightarrow  P
\stackrel{\beta}\rightarrow  M  \rightarrow 0$$ we obtain the
short exact sequence:
$$(\ast\ast)\quad 0\rightarrow  Hom(M,Q)  \stackrel{\circ\beta}\rightarrow
Hom(P,Q) \stackrel{\circ\alpha}\rightarrow  Hom(M,Q)  \rightarrow
0$$ since $Ext(M,Q)=0$ (\cite[Proposition 2.3]{Holm}). On the
other hand, $u\in Hom(M,Q)$. Then, from the exactness of
$(\ast\ast)$, there is a morphism $\upsilon: P\rightarrow Q$ such
that $\upsilon\circ \alpha=u$. Consequently, we can verify that
the morphism $\gamma:P\rightarrow Q\oplus Q$ defined by
$\gamma(p):=(\upsilon(p),u\circ\beta(p))$ whenever $p\in P$ is the
desired morphism.
\end{proof}

Dually, we obtain easily the injective version of Lemma
\ref{lemma1}.\\

In section 2 and 3, motived by the important role of the rings of
global and weak dimensions smaller or equal to one in several
areas of algebra, we study the rings of Gorenstein homological
dimensions smaller or equal to one which we call, by analogy to
the classical ones, Gorenstein hereditary, semihereditary rings,
specially particular cases of these rings which we call strongly
Gorenstein
hereditary, semihereditary rings.\\

\section{On (strongly) Gorenstein hereditary rings}

The aim of this section is to characterize the rings of Gorenstein
global dimension smaller or equal than 1, specially a particular
case of them over which every Gorenstein projective module is
strongly Gorenstein projective.

\begin{definitions}\
\begin{enumerate}
    \item A ring $R$ is called  a Gorenstein hereditary ring ($G$-hereditary for short)
    if every submodule of a projective module is $G$-projective (i.e., $G-gldim(R) \leq 1$) and  $R$ is called a Gorenstein Dedekind ring ($G$-Dedekind for short), if it is  a $G$-hereditary domain.
      \item A ring $R$ is called a strongly Gorenstein hereditary ring
      ($SG$-hereditary for short)
      if every submodule of a projective module is
      $SG$-projective and $R$ is called strongly Gorenstein Dedekind ($SG$-Dedekind for short), if it is an $SG$-hereditary
domain.
\end{enumerate}
\end{definitions}
\begin{remark}
It is easy to see that the definition of a strongly Gorenstein
hereditary ring is equivalent to say that every submodule of a
strongly Gorenstein projective module is strongly Gorenstein
projective.
\end{remark}
In the next, we give a characterization of the G-hereditary rings.

\begin{proposition}\label{G-hereditary} Let $R$ be a ring with finite Gorenstein global dimension. The following assertions are equivalent:
\begin{enumerate}
    \item $R$ is G-hereditary.
    \item $Gpd_R(M)\leq 1$ for all finitely generated R-modules M.
    \item Every ideal of $R$ is Gorenstein
    projective.
    \item $id_R(P)\leq  1$ for all $R$-modules $P$ with finite  $pd_R(P)$.
    \item $id_R(P')\leq  1$ for all projective $R$-modules $P'$.
    \item $pd_R(E)\leq  1$ for all $R$-modules $E$ with finite  $id_R(E)$.
    \item $pd(E') \leq 1$ for all injective $R$-modules $E'$.
\end{enumerate}
\end{proposition}
\begin{proof} All no obvious implications follow immediately from
\cite[Theorem 1.1]{Bennis and Mahdou2}, \cite[Theorems 2.20 and
2.22]{Holm} and \cite[Lemma 9.11]{Rotman}.
\end{proof}

The main result, in this section, is the following
characterization of
 the $SG$-hereditary rings.
 \begin{theorem}\label{caracterisation1}
Let $R$ be a ring. The following assertions are equivalent:
\begin{enumerate}
    \item $R$ is $SG$-hereditary.
    \item For every $R$-module $M$, there exists a short exact
    sequence: $$0\longrightarrow M \longrightarrow Q \longrightarrow M \longrightarrow 0$$ where $pd_R(Q)\leq 1$, and for every projective module
    $P$ there is an integer $i> 1$ such that  $Ext_R^i(M,P)=0$.
     \item For every $R$-module $M$, there exists a short exact
    sequence: $$0\longrightarrow M \longrightarrow E \longrightarrow M \longrightarrow 0$$ where $id_R(E)\leq 1$, and for every injective module
    $I$ there is an integer $i> 1$ such that  $Ext_R^i(I,M)=0$.
\end{enumerate}
\end{theorem}

\begin{proof}
$1\Rightarrow 2.$ Assume $(1)$ and we claim $(2)$. Let $M$ be an
arbitrary $R$-module. Pick a short exact sequence $0\rightarrow G
\rightarrow Q \rightarrow M \rightarrow 0$ where $Q$ is a
projective module. The module $G$ is immediately strongly
Gorenstein projective by the hypothesis conditions. Hence, from
\cite[Proposition 2.9]{Bennis and Mahdou1}, there exists  a short
exact sequence $0\rightarrow G \rightarrow P \rightarrow
G\rightarrow 0$ where $P$ is projective. Consider the following
diagram:
$$\begin{array}{ccccccc}
  0\rightarrow & G & \rightarrow & P & \rightarrow & G & \rightarrow 0 \\
    &  \downarrow & &  && \downarrow &  \\
  0\rightarrow & Q & \stackrel{\iota}\rightarrow & Q\oplus Q & \stackrel{j}\rightarrow & Q & \rightarrow 0 \\
\end{array}$$
By Lemma \ref{lemma1}, the above diagram  can be completed and so
applying the Snake Lemma and the fact that $M\cong coker
(G\rightarrow Q)$, we can construct a short exact sequence
$0\rightarrow M \rightarrow X \rightarrow M \rightarrow 0$ where
$X\cong coker(P\rightarrow Q\oplus Q)$. Clearly, $pd(X)\leq 1$.
Thus,  we have the desired short exact sequence. Furthermore,
since $G.gldim(R)\leq 1$, for any module $M$ and any projective
module $P$ we have $Ext^i(M,P)=0$ for any integer $i>1$.

$2\Rightarrow 1.$ Assume the second assertion and let $M$ be a
submodule of a projective $R$-module $P$. We claim that $M$ is
strongly Gorenstein projective. Applying the hypothesis conditions
to the module $P/M$, there exists a short exact sequence
$$0\rightarrow P/M\rightarrow Q\rightarrow P/M
\rightarrow 0$$ where $pd_R(Q)\leq 1$. Now, consider the following
diagram with exact rows and columns:
\begin{center}
$\begin{array}{ccccccc}
   & 0 &  & 0 &  & 0 &  \\
   &\downarrow & & \downarrow & & \downarrow &  \\
  0\rightarrow & M & \rightarrow &X & \rightarrow & M & \rightarrow 0 \\
  & \downarrow &  & \downarrow &  & \downarrow &  \\
0\rightarrow &  P & \rightarrow & P\oplus P & \rightarrow & P & \rightarrow 0 \\
   & \downarrow &  & \downarrow &  & \downarrow &  \\
  0\rightarrow &  P/M& \rightarrow & Q & \rightarrow & P/M & \rightarrow 0 \\
  &\downarrow &  & \downarrow &  & \downarrow &  \\
   & 0 &  & 0 &  & 0 &  \\
\end{array}$\end{center}
From the middle vertical short exact sequence, we deduce that $X$
is projective. Moreover, by the hypothesis conditions again
(applying to the module $M$), for every projective module $F$
there is  an integer $i> 1$ such that $Ext_R^i(M,F)=0$. Then, from
the short exact sequence $0\rightarrow M \rightarrow X \rightarrow
M \rightarrow 0$ we get  $Ext_R(M,F)=0$. So, from
\cite[Proposition 2.9]{Bennis and Mahdou1}, $M$ is an
$SG$-projective module, as desired.\\

$1\Rightarrow 3.$ Using the dual of the results in the proof of
the implication $1\Rightarrow 2$, the proof of the present
implication is similar to the one of $1\Rightarrow 2$.

$3\Rightarrow 1.$ Assume $(3)$. By a dual argument to the one of
the first part of the implication $2\Rightarrow 1$ we can prove
that $I/M$ is a strongly Gorenstein injective module for an
arbitrary module $M$ and an injective module $I$ which contains
$M$. Hence, $Gid(M)\leq 1$. Consequently, by \cite[Theorem
1.1]{Bennis and Mahdou2}, $G.gldim(R)=GID(R)\leq 1$, and so, $R$
is $G$-hereditary. Let $M$ be a submodule of a projective module.
We claim  that $M$ is strongly Gorenstein projective. Hence, $M$
is Gorenstein projective since $R$ is a Gorenstein hereditary
ring. Thus, from \cite[Theorem 2.20]{Holm}, $Ext(M,P)=0$ for any
projective module $P$. Moreover, by the hypothesis conditions,
there is a short exact sequence $(\star)\quad 0\rightarrow M
\rightarrow E \rightarrow M \rightarrow 0$ where $id(E)\leq 1$.
Hence, $pd(E)\leq 1$ (by Proposition \ref{G-hereditary}). Applying
\cite[Theorem 2.5]{Holm} to $(\star)$, we conclude that $E$ is
Gorenstein projective. Thus, from \cite[Propositon 2.27]{Holm},
$E$ is projective. Consequently, from \cite[Proposition
2.9]{Bennis and Mahdou1}, $M$ is strongly Gorenstein projective,
as desired.
\end{proof}
\begin{corollary}\label{cor1}
Let $R$ be a ring with finite Gorenstein global dimension. The
following assertions are equivalent :
\begin{enumerate}
    \item $R$ is $SG$-hereditary.
    \item For every $R$-module $M$, there exists a short exact sequence  $$0
    \longrightarrow M \longrightarrow Q \longrightarrow M
    \longrightarrow 0$$ such that  $pd_R(Q)\leq1$.
   \item For every $R$-module $M$, there exists a short exact sequence  $$0
    \longrightarrow M \longrightarrow E \longrightarrow M
    \longrightarrow 0$$ such that  $id_R(E)\leq1$.
\end{enumerate}
\end{corollary}
\begin{proof} From \cite[Theorems 2.20 and 2.22]{Holm} and the
definition of $G.gldim(-)$ (see \cite[page 1]{Bennis and
Mahdou2}), for every $R$-module and every projective $R$-module
$P$ we have $Ext_R^n(M,P)=0$ where $n=G.gldim(R)$; and similarly
for every injective $R$-module $I$ we have $Ext_R^n(I,M)=0$. Then,
this Corollary  follows directly from Theorem
\ref{caracterisation1}.

\end{proof}

It's  well-known that the hereditary rings (resp. Dedekind
domains) are coherent (resp. Noetherian). Now, it is natural to
ask  what about the $G$-hereditary rings. Now, we can give an
affirmative answer just in the strongly Gorenstein hereditary
case.

\begin{theorem}\label{th coherent}\
\begin{enumerate}
    \item Every $SG$-hereditary ring is coherent.
    \item Every coherent $G$-Dedekind domain, in particular  every $SG$-Dedekind domain,  is Noetherian.
\end{enumerate}
\end{theorem}

\begin{proof}[Proof of Theorem \ref{th coherent}]
 $\mathbf{(1)}$ Assume that $R$ is an $SG$-hereditary ring and let $I$ be a finitely
generated ideal of $R$. Then, $I$ is an $SG$-projective
$R$-module, since $I$ is a submodule of the projective $R$-module
$R$. Therefore, $I$ is a finitely presented $R$-module (by
    \cite[Theorem 3.9]{Bennis and Mahdou1}). So, $R$ is coherent, as desired.\\

$\mathbf{(2)}$ Let $R$ be a $G$-Dedekind coherent domain. Then,
$G-gldim(R)\leq 1$. If $G-gldim(R)= 0$, then $R$ is a
quasi-Frobenius ring (by \cite[Proposition 2.6]{Bennis and
Mahdou2}) and so is Noetherian. Now, suppose that $G-gldim(R)=1$.
The finitistic Gorenstein projective dimension of $R$, denoted by
$FGPD(R)$ is finite (See \cite[p.182]{Holm}). Namely, $FGPD(R) =
G-gldim(R)$. From \cite[Theorem 2.28]{Holm}, it is equal to the
finitistic projective dimension of $R$, denoted by $FPD(R)$. Then,
$FPD(R)=1$. Therefore, $R$ is Noetherian (from \cite[Theorems
2.5.14]{Glaz}).\\

Taking the consideration $(1)$ above, the particular case (i.e;
where $R$ is $SG$-Dedekind domain) is immediate.
\end{proof}

\begin{proposition}\label{G-flat SG-hereditary}
Let $R$ be a Gorenstein hereditary ring. The following statements
are equivalents:
\begin{enumerate}
    \item $R$ is a strongly Gorenstein hereditary ring;
    \item Every $G$-projective $R$-module is strongly Gorenstein
    projective;
    \item $R$ is coherent and every Gorenstein flat $R$-module is strongly Gorenstein
    flat module.
\end{enumerate}
\end{proposition}
\begin{proof}
$1\Leftrightarrow 2.$ Obvious since every $G$-projective module is
a submodule of a projective module.\\
$2\Rightarrow 3.$  The coherence of $R$ is guarantied by Theorem
\ref{th coherent}(1). Now, let $M$ be a $G$-flat module. By
hypothesis $ Gpd_R(M)\leq 1$. Thus, from Theorem
\ref{caracterisation1}, there is an exact sequence $0\rightarrow M
\rightarrow X \rightarrow M \rightarrow 0$ where $pd(X)\leq 1$.
Then, $id(Hom_\mathbb{Z}(X,\mathbb{Q}/\mathbb{Z}))=fd(X)\leq
pd(X)\leq 1$. Furthermore, from \cite[Theorem 3.7]{Holm}, $X$ is a
Gorenstein flat module since $\mathcal{GF}(R)$ is projectively
resolving. Hence, by \cite[Proposition 3.11]{Holm},
$Hom_\mathbb{Z}(X,\mathbb{Q}/\mathbb{Z})$ is Gorenstein injective.
Consequently, by the dual of \cite[Proposition 2.27]{Holm},
$Hom_\mathbb{Z}(X,\mathbb{Q}/\mathbb{Z})$ is injective. Then, from
\cite[Theorem 1.2.1]{Glaz}), $X$ is flat. Hence,  $M$ is
immediately $SG$-flat (by \cite[Proposition 3.6]{Bennis and
Mahdou1} and since for any injective module $I$, we have
$Tor(M,I)=0$ as $M$ is
$G$-flat).\\
$3\Rightarrow 2$. Let $I$ be an injective $R$-module. From
\cite[Corollary 2.7]{Bennis and Mahdou2}, $fd_R(I)\leq1$. Then,
from \cite[Theorem 3.8]{Ding}, $R$ is an $1-FC$ ring (i.e.,
coherent ring with $Ext^2_R(P,R)=0$ for each finitely presented
$R$-module $P$). Now,  let $M$ be a $G$-projective $R$-module.
Then, $M$ embeds in a projective $R$-module. So, from
\cite[Theorem 7]{Chen}, $M$ is $G$-flat. Then, by hypothesis $M$
becomes $SG$-flat. Hence, there exists a short exact sequence
$0\rightarrow M \longrightarrow F \rightarrow M \rightarrow 0$
where $F$ is flat. By the resolving of
 the class $\mathcal{GP}(R)$ and from the short exact sequence above we
deduce that  $F$ is $G$-projective (since $M$ is $G$-projective).
On the other hand, $pd_R(F)<\infty$ (by \cite[Corollary
2.7]{Bennis and Mahdou2} and since $F$ is flat). Therefore, $F$ is
projective by \cite[Proposition 2.27]{Holm}. So $M$ is
$SG$-projective (by \cite[Proposition 2.9]{Bennis and Mahdou1} and
since $Ext(M,P)=0$ for every projective module $P$ as  $M$ is
$G$-projective).
\end{proof}

\begin{remark}
Using \cite[Corollary 1.2(2)]{Bennis and Mahdou2}, we say clearly
that a $G$-hereditary ring (in particular an $SG$-hereditary ring)
is hereditary if, and only if, $wdim(R)$ is
    finite.
\end{remark}

In what follows we give an example of non $G$-semisimple
$SG$-hereditary ring and  a $G$-hereditary ring which is not
$SG$-hereditary.
\begin{example}
Consider a non-semisimple quasi-Frobenius rings $R = K[X]/(X^2$),
$R' = K[X]/(X^3)$ where $K$ is a field, and a non- Noetherian
hereditary ring $S$. Then,
\begin{enumerate}
    \item $R\times S$ is a strongly Gorenstein hereditary ring which is
    not hereditary.
    \item $R'\times S$ is a Gorenstein hereditary ring which is
    not strongly Gorenstein hereditary.
\end{enumerate}
\end{example}
\begin{proof}
From [4, Example 3.4], the rings $R\times S$ and $R'\times S$ are
 both Gorenstein hereditary rings with infinite weak dimension.\\
$(1)$ We have to prove that $R\times S$ is strongly Gorenstein
hereditary. From Proposition \ref{G-flat SG-hereditary}, it
remains to prove that every Gorenstein projective module is
strongly Gorenstein projective. Let $M$ be a Gorenstein projective
$R\times S$-module. We claim that $M$ is an $SG$-projective
module. We have the isomorphism of $R\times S$-modules:
$$M\cong M\otimes_{R\times
    S}R\times S\cong M\otimes_{R\times
    S}(R\times0\oplus 0\times S)\cong M_1\times M_2$$ where
    $M_1=M\otimes_{R\times
    S} R$  and  $M_2=M\otimes_{R\times
    S} S$ (for more details see \cite[p.102]{Berrick}). By \cite[Lemma 3.2]{Bennis and Mahdou3},   $M_1$ (resp. $M_2$) is
    a $G$-projective $R$-module (resp. $S$-module). Then, since $R$ is strongly Gorenstein semisimple and $S$ is hereditary, $M_1$ (resp. $M_2$) is an $SG$-projective $R$-module (resp. $S$-module) (precisely $M_2$ is a projective $S$-module).
     On the other hand, the family $\{R, S\}$ of rings satisfies the conditions of
    \cite[Lemma 3.3]{Bennis and Mahdou3} (by \cite[Corollary 2.7]{Bennis and Mahdou2}  since $G.gldim(R)$ and $G.gldim(S)=gldim(S)$ are  finite). Thus,   $M=M_1\times M_2$ is
   an  $SG$-projective  $R\times S$-module, as desired.\\
   $(2)$ We have to prove that $R'\times S$ is not strongly Gorenstein
   module. By \cite[Corollary 3.10]{Ouarghi}, there exists a Gorenstein projective $R'$-module $M$ which is not strongly Gorenstein projective.
   And by, \cite[Lemma 3.2]{Ouarghi}, $M\times S$ is Gorenstein
projective $R'\times
   S$-module which is not strongly Gorenstein projective. Thus,
   from Proposition \ref{G-flat SG-hereditary}, $R'\times S$ is
   not strongly Gorenstein hereditary.
   \end{proof}
\section{On (strongly) Gorenstein semihereditary rings}

The aim of this section is to characterize the rings of Gorenstein
weak  dimension smaller or equal than 1, specially a particular
case of them over which every Gorenstein flat module is strongly
Gorenstein flat.

\begin{definition}\
\begin{enumerate}
    \item A ring $R$ is called Gorenstein semihereditary ($G$-semihereditary fort short) if $R$ is coherent
and every submodule of flat module is $G$-flat (i.e., $R$ is
coherent and $G.wdim(R) \leq 1$).
    \item A ring $R$ is called strongly Gorenstein semihereditary ($SG$-semihereditary for short) if $R$ is coherent
and every submodule of flat module is $SG$-flat.
\end{enumerate}

\end{definition}
\begin{remarks}\label{semiher resol}$(i)$ Clearly we have the followings equivalences:
\begin{enumerate}
    \item A $G$-semihereditary ring $R$ is semihereditary if, and
only if, $G.wdim(R)$ is finite (\cite[Corollary 1.2(3)]{Bennis and
Mahdou2}).
    \item A $G$-semihereditary ring $R$ is $SG$-semihereditary if, and
only if, every $G$-flat module is $SG$-flat.
\end{enumerate}
$(ii)$ Every $SG$-hereditary ring is $SG$-semihereditary (by
    \cite[Corollary 1.2(1)]{Bennis and Mahdou2}, Proposition \ref{G-flat
    SG-hereditary} and $(i2)$ above).\\
$(iii)$ Every Noetherian $G$-semihereditary (resp.
$SG$-semihereditary) ring is $G$-hereditary (resp.
$SG$-hereditary) (From \cite[(Theorem 12.3.1]{Enochs3} and also
Proposition \ref{G-flat
    SG-hereditary} and $(i2)$ above in the strongly case).

\end{remarks}
Recall that we say that  an $R$-module M has $FP$-injective
dimension at most $n$ (for some $n \geq 0$) over a ring $R$,
denoted by $FP-id_R(M)\leq n$, if  $Ext_R^{n+1} (P,M) = 0$ for all
finitely presented $R$-modules $P$. Recall also that  $R$ is
called $n-FC$ (for some $n \geq 0$), if it is coherent and it has
self-$FP$-injective dimension at most $n$ (i.e., $FP-id_R(R) \leq
n$). Now, we give a characterization of the  $G$-semihereditary
rings.
\begin{proposition}\label{caracterisation de G-semiherdiatry} Let R be a coherent ring, then the following statements are equivalent:
\begin{enumerate}
    \item R is $G$-semihereditary.
    \item $Gfd_R(M) \leq 1$ for all finitely presented R-modules M.
    \item $Gpd_R(M) \leq 1$ for all finitely presented R-modules M.
    \item Eery finitely generated ideal of $R$ is Gorenstein flat.
    \item $fd_R(I) \leq1$ for all injective $R$-modules I.
    \item $FP-id_R(F) \leq 1$  for all flat $R$-modules F.
    \item $fd_R(E) \leq 1$  for all $FP$-injective $R$-modules E.
\end{enumerate}
\end{proposition}
\begin{proof} First, note that a  ring $R$ is $1-FC$ if, and only if,
$R$ is
$G$-semihereditary (by \cite[Theorem 7]{Chen}).\\
The implications $(1)\Rightarrow (2) \Rightarrow (4)$ are obvious.\\
$(4)\Rightarrow (5).$ Follows from \cite[Theorem 3.14]{Holm} and \cite[Theorem 1.3.8]{Glaz}.\\
$(5)\Rightarrow (6)\Rightarrow (7) \Rightarrow (1).$ Follows from \cite[Theorem 3.8]{Ding}.\\
$(1)\Rightarrow(3).$ Follows from \cite[Theorem 7]{Chen}.
\end{proof}

Now, we give the main two results in this section.

\begin{theorem}\label{thm semihere1}
Let $R$ be a ring such that every direct limit of $SG$-flat
$R$-modules is $SG$-flat. Then, the following statements are
equivalent:
\begin{enumerate}
    \item $R$ is $SG$-semihereditary.
    \item Every finitely generated submodule  of a projective
    module is $SG$-projective.
\end{enumerate}
\end{theorem}

\begin{proof}
 We assume that $R$ is an $SG$-semihereditary ring and
let $M$ be a finitely generated submodule  of a projective
    module $P$. Then, $M$ is a finitely presented $SG$-flat module
    (since $R$ is $SG$-semihereditary and $R$ is coherent). Then,
    $M$ is $SG$-projective (By
    \cite[Proposition 3.9]{Bennis and Mahdou1}), as desired.\\
 Conversely, we assume that every finitely generated
submodule  of a projective module is $SG$-projective. Our aim is
to show that $R$ is $SG$-semihereditary. Let $I$ be a finitely
generated ideal of $R$. By hypothesis, $I$ is a finitely generated
$SG$-projective $R$-module. Then, by \cite[Proposition 3.9]{Bennis
and Mahdou1}, $I$ is a finitely presented $SG$-flat $R$-module.
So, by Proposition \ref{caracterisation de G-semiherdiatry},  $R$
is a
 $G$-semihereditary ring (see that $R$ is coherent). Now, to prove
that $R$ is $SG$-semihereditary ring, it suffices, by Remark
\ref{semiher resol}, to prove that every $G$-flat $R$-module is
$SG$-flat. So, let $M$ be a $G$-flat $R$-module. Then, $M$ embeds
in a flat $R$-module $F$. By Lazard's theorem (\cite[Theorem
1.2.6]{Glaz}), there is a direct system $(L_i,\varphi_{i,j}
)_{i\in I}$ of finitely generated free $R$-modules such that
$\underrightarrow{lim} (L_i,\varphi_{i,j})\cong F$. On the other
hand $\underrightarrow{lim} (L_i,\varphi_{i,j})= \displaystyle
\frac{\oplus L_i}{S}$ where $S$ is the submodule generated by all
elements $\lambda _j\circ \varphi _{i,j}(a_i)-\lambda _i (a_i)$,
where $a_i\in L_i$ and $i\leq j$, and for each $i\in I$ the
homomorphism $\lambda _i$ is the injection of $L_i$ into the sum
$\oplus
L_i$ (for more details see \cite[pages 32, 33 and 34]{Atiyah}).\\
We can identify $M$ to a submodule of $\displaystyle\frac{\oplus
L_i}{S}$ and we consider an $R$-module $A$ and an homomorphism
$\alpha$ of $R$-modules such that the short sequence of
$R$-modules $$0\longrightarrow M \hookrightarrow \displaystyle
\frac{\oplus L_i}{S} \stackrel{\alpha} \longrightarrow A
\longrightarrow 0$$ is exact.\\Now, consider the family of exacts
sequences $0\longrightarrow M_i \hookrightarrow \displaystyle L_i
\stackrel{\alpha \circ \bar{\lambda_i}} \longrightarrow
A_i\longrightarrow 0$ where $M_i=ker(\alpha \circ
\bar{\lambda_i})$, $A_i=Im (\alpha \circ \bar{\lambda_i})$ and the
homomorphism $\overline{\lambda_i}: L_i \mapsto \displaystyle
\frac{\oplus L_i}{S}$ is such that for each $a\in L_i$,
$\overline{\lambda_i} (x)=\overline{\lambda_i(x)}$.\\
$\mathbf{a)}$ First, we claim that $A=\underrightarrow{lim}
    (A_i, \subseteq)$. For each  $x\in A_i$, there exists an element $y\in L_i$
such that $\alpha \circ \bar{\lambda_i}(y)=x$. By definition of
direct system, we have $\overline{\lambda_i}=\overline{\lambda_j}
\circ \varphi_{i,j}$ (for $i\leq j$; $i,j\in I$). Thus, we deduce
that
$$x=\alpha \circ \bar{\lambda_i}(y)=\alpha \circ
\bar{\lambda_j}(\varphi_{i,j}(y)) \in \alpha \circ
\bar{\lambda_j}(L_j)=A_j$$ Consequently, for
$i\leq j$, we have $A_i\subseteq A_j$. So, we conclude that\\
$\underrightarrow{lim}(A_i,\subseteq)=\sum A_i$. On the other
hand, for every $i\in I$,
$$A_i=\alpha\circ\overline{\lambda_i}(L_i)=\alpha(\overline{\lambda_i}(L_i))\subseteq
\alpha\displaystyle(\frac{\oplus L_i}{S})=A$$ This implies that
$\underrightarrow{lim} (A_i)=\sum A_i
    \subseteq A$.\\ Conversely, for each $x\in A$, there
    exists  $y\in \displaystyle \frac{\oplus L_i}{S}$ such that
    $\alpha(y)=x$. We have  $y=\overline{(x_i)_{i\in I}}=\overline{\sum
    \lambda_i(x_i)}=\sum\overline{
    \lambda_i(x_i)}$ such that $x_i=0$ except for a finite elements of $I$. Then, $x=\alpha(y)=\sum\alpha\circ
    \overline{\lambda_i}(x_i) \in \sum A_i.$ Thus, we conclude
    that:
    $A=\underrightarrow{lim}
    (A_i)$.\\

$\mathbf{b)}$ For each $x\in M_i=Ker(\alpha \circ
\overline{\lambda_i})$, we have $\alpha \circ
\overline{\lambda_j}(\varphi_{i,j}(x))=\alpha \circ
\overline{\lambda_i}(x)=0$. So, $\varphi_{i,j}(x)\in M_j$. Then,
for $i\leq j$, the family of homomorphisms: $\varphi_{i,j}': M_i
\mapsto M_j$, such that for
    each $x\in M_i\subseteq L_i$,
    $\varphi_{i,j}'(x)=\varphi_{i,j}(x)$ are well defined and the system $(M_i,\varphi_{i,j}')_{i\in I}$ is
    direct. More thus,
    the following diagram:
   $$ (\star)\begin{array}{ccccccc}
      0 \longrightarrow & M_{i} & \hookrightarrow & L_{i}  &\stackrel{\alpha \circ \overline{\lambda_{i}}} \longrightarrow & A_{i} &\longrightarrow 0 \\
      & \varphi_{i,j}'\downarrow &  &  \varphi_{i,j}\downarrow &  &\mu_{i,j} \downarrow  &  \\
      0\longrightarrow & M_{j} &\hookrightarrow &L_{j}  &\stackrel{\alpha \circ \overline{\lambda_{j}}}\longrightarrow & A_{j} &\longrightarrow 0 \\
    \end{array}$$ where $\mu_{i,j}$ is the embedding of $A_i$ in $A_j$, is commutative. So, the short sequence of direct
    system over $I$ induced from $(\star)$ :$$0\longrightarrow
    (M_i,\varphi_{i,j}')_{i\in I}\longrightarrow (L_i,\varphi_{i,j})_{i\in I}
    \longrightarrow(A_i,\subseteq )_{i\in I} \longrightarrow 0$$ is
exact and by \cite[Exercice 18, p.33]{Atiyah},
$\underrightarrow{lim}(\alpha \circ
\overline{\lambda_{i}})=\alpha$.
    Consequently, the short exact sequence:$$ 0 \longrightarrow \underrightarrow{lim}
    (M_i,\varphi_{i,j}')\longrightarrow \displaystyle \frac{\oplus L_i}{S} \stackrel{\alpha}\longrightarrow A\longrightarrow 0.$$\\
 is exact and so   $\underrightarrow{lim}
    (M_i)\cong ker (\alpha) =M$.\\
$\mathbf{c)}$ For each $i\in I$, $M_i$ is direct limit of his
finitely generated
    submodules $(M_i^j)_j$, and for each $j$, $M_i^j\subseteq M_i \subseteq L_i$. Then, by hypothesis, $M_i^j$ is $SG$-projective. So, $M_i$ is direct limit of a finitely generated   $SG$-projective modules (then $SG$-flat modules by \cite[Proposition 3.9]{Bennis and Mahdou1}). Then, by hypothesis,  $M_i$ is
    $SG$-flat.\\
$\mathbf{conclusion:}$ By (b) and (c) we conclude that $M$ is a
direct limit of $SG$-flat
     modules. Thus, by hypothesis, $M$ is
    $SG$-flat, as desired.

\end{proof}
\begin{theorem}
Let $R$ be a ring such that every $G$-flat module is $SG$-flat.
Then, $R$ is $SG$-semihereditary if, and only if, every finitely
generated ideal is $SG$-projective.
\end{theorem}
\begin{proof} We assume that $R$ is $SG$-semihereditary and let $I$
be a finitely generated ideal of $R$. Then, $I$ is a finitely
presented $SG$-flat $R$-module (since $R$ is coherent). Thus, by
\cite[Proposition 3.9]{Bennis and Mahdou1}, $I$ is $SG$-projective.\\

Conversely, we assume that every finitely generated ideal of $R$
is $SG$-projective and every $G$-flat module is $SG$-flat. It is
clear that R is coherent (by \cite[Proposition 3.9]{Bennis and
Mahdou1}, we deduce that every finitely generated ideal is
finitely presented G-flat). We have also that $Gfd(R/I)\leq 1$ for
every finitely generated ideal $I$. Then, for every injective
module E we have $Tor^2(E,R/I) = 0$ (by \cite[Theorem
3.14]{Holm}). Therefore, by \cite[Theorem 1.3.8]{Glaz}, $fd(E)\leq
1$. Using \cite[Theorem 3.8]{Ding}, $FP -id(R) \leq 1$ and so $R$
is $1-FC$ since it is coherent. Then, by \cite[Theorem 7]{Chen},
$Gfd(M) \leq 1$ for every module $M$ and then $G.wdim(R) \leq 1$.
Moreover, every Gorenstein flat module is strongly Gorenstein
flat. Then R is $SG$-semihereditary (by Remarks \ref{semiher
resol}(i2)).
\end{proof}
%%%%%%%%%%%%%%%%%%%%%%%%%%%%%%%%%%%%%%%%%%%%%%%%%%%%%%%%%
%%%%%%%%%%%%%%%%%%%%%%%%%%%%%%%%%%%%%%%%%%%%%%%%%%%%%%%%%
%%%REFERENCES%%%%%%%%%%%%%%%%%%%%%%%%%%%%%%%%%%%%%%%%%%%%
%%%%%%%%%%%%%%%%%%%%%%%%%%%%%%%%%%%%%%%%%%%%%%%%%%%%%%%%

%%%%%%%%%%%%%%%%%%%%%%%%%%%%%%%%%%%%%%%%%%%%%%%%%%%%

%%%%%%%%%%%%%%%%%%%%%%%%%%%%%%%%%%%%%%%%%%%%%%%%%%%%

\end{document}